\documentclass{mcom-l}

\newcommand{\R}{\mathbb{R}}
 
\newcommand{\Z}{\mathbb{Z}}
\newcommand{\N}{\mathbb{N}}

\usepackage{amssymb}
\usepackage{hyperref}
\usepackage{cleveref}
\usepackage{tikz}
\usetikzlibrary{arrows,backgrounds}
\usepackage{algpseudocode}
\algtext*{EndWhile}
\algtext*{EndFor}
\algtext*{EndIf}

\copyrightinfo{}{}

\newtheorem{theorem}{Theorem}[section]
\newtheorem{prop}[theorem]{Proposition}
\newtheorem{lemma}[theorem]{Lemma}
\newtheorem{corollary}[theorem]{Corollary}
\newtheorem{algorithm}[theorem]{Algorithm}

\theoremstyle{definition}
\newtheorem{definition}[theorem]{Definition}
\newtheorem{example}[theorem]{Example}

\newtheorem{remark}[theorem]{Remark}

\numberwithin{equation}{section}

\begin{document}

\title[Factorization as Subset-Sum]{Integer Factorization as Subset-Sum Problem}


\author[M. Hittmeir]{Markus Hittmeir}
\address{}
\curraddr{}
\email{mhittmeir@sba-research.org}
\thanks{SBA Research (SBA-K1) is a COMET Centre within the framework of COMET – Competence Centers for Excellent Technologies Programme and funded by BMK, BMDW, and the federal state of Vienna. The COMET Programme is managed by FFG}


\date{}

\dedicatory{}

\begin{abstract}
This paper elaborates on a sieving technique that has first been applied in 2018 for improving bounds on deterministic integer factorization. We will generalize the sieve in order to obtain a polynomial-time reduction from integer factorization to a specific instance of the multiple choice subset-sum problem. As an application, we will improve upon special purpose factorization algorithms for integers composed of divisors with small difference. In particular, we will refine the runtime complexity of Fermat's factorization algorithm by a large subexponential factor. Our first procedure is deterministic, rigorous, easy to implement and has negligible space complexity. Our second procedure is heuristically faster than the first, but has non-negligible space complexity. 
\end{abstract}

\maketitle

\section{Introduction}
Let $N$ be a natural number. We consider the problem to compute the prime factorization of $N$. No efficient general-purpose algorithm for integer factorization is known. The assumed computational hardness of this problem is used in cryptographic applications, e.g. in the RSA-cryptosystem. The fastest currently known algorithms for computing the prime factorization of natural numbers achieve sub-exponential runtime complexities. In practice, methods such as the Number Field Sieve, the Quadratic Sieve and the Elliptic Curve Method are used for factoring large integers. The monographs \cite{Rie} and \cite{Wag} provide an overview on the variety of techniques. However, these methods are either conditional in some sense or rely on probabilistic and heuristic arguments. Among such approaches, Shor's algorithm (\cite{sh}) takes in an exceptional position. It runs in polynomial-time, but relies on quantum computing. 

Another class of algorithms concerns deterministic factorization with rigorous runtime complexity bounds. All known bounds in this area are still exponential in the length of the input number. In 1974, Lehman (\cite{Leh}) published an algorithm that runs in $\widetilde{O}(N^{1/3})$ in big O tilde notation, which ignores logarithmic factors in the runtime bound. Compared to the bound $\widetilde{O}(N^{1/2})$ of the trial division algorithm, this was a major improvement. A few years later, Pollard and Strassen used fast polynomial multiplication and multipoint evaluation techniques to establish a deterministic and rigorous factorization algorithm running in time $\widetilde{O}(N^{1/4})$. Since the publication of Strassen's paper in 1977, there have been a few refinements of the runtime bound, the best of which has been proved in 2014 (\cite{CosHar}). 
In 2018, the first superpolynomial improvement has been achieved via a new sieve technique (\cite{Hit1}). The runtime of the resulting algorithm is given by 
\[
\widetilde{O}\left(N^{1/4}\exp(-C\log N/\log\log N)\right)
\] 
for a positive constant $C$. Then, in 2020 (\cite{Hit2}), we finally improved the $1/4$-exponent threshold by achieving a complexity of $\widetilde{O}(N^{2/9})$ for deterministic factorization. The main idea of this algorithm is based on combining fast polynomial arithmetic techniques with Lehman's method mentioned above. Shortly thereafter, Harvey (\cite{Har}) showed that our algorithm can be modified to obtain the complexity $\widetilde{O}(N^{1/5})$. The most recent refinement (\cite{HarHit}) has been proved in 2021.

To summarize, there appear to be two different strands of research on the integer factorization problem:
\begin{enumerate}
	\item{The study of algorithms intended for practice: Here, we are more interested in usability than in generality or rigorosity of our proofs. In addition to the already mentioned methods with sub-exponential complexity, there are also special purpose methods that work well for numbers with certain properties (e.g., Fermat's factorization algorithm or Pollard's $P-1$-algorithm)}
	\item{The study of integer factorization as a theoretical problem: Besides its crytographical significance, integer factorization is a central problem in algorithmic number theory. In this context, we try to find the most efficient way to factorize any natural number on a deterministic Turing machine. The runtime analysis should be a rigorous mathematical proof.}
\end{enumerate}

The general goal of the present paper is to elaborate on the techniques from the recent improvements in the second area and prepare them for application in the first area. In particular, we will consider the sieve technique established in \cite{Hit1}, which we call the \emph{hyperbolic sieve}. 
Let us discuss this idea in greater detail. We consider the case $N=pq$, where $p$ and $q$ are distinct primes and $p<q$. Clearly, $p$ and $q$ are solutions to the congruential equation $N\equiv xy \mod m$ for every $m\in\N$ which is coprime to $N$. Throughout this paper, let $\Z_n:=\Z/n\Z$ for every $n\in\N$. We then define the set 
\[
\mathcal{H}_{N,m}:=\{(x,y)\in\Z_m: N\equiv xy\mod m\}
\]
and note that $(p\mod m, q\mod m)\in \mathcal{H}_{N,m}$. These sets are called \emph{modular hyperbolas} and have been studied extensively in \cite{Shp}. It is easy to see that the number of pairs in $\mathcal{H}_{N,m}$ equals $\phi(m)$, where $\phi$ is Euler's totient function. In order to obtain applicable information about $p$ and $q$, we consider the corresponding set
$
\mathcal{L}_{N,m}:=\{x+y\pmod m: (x,y)\in\mathcal{H}_{N,m}\}.
$
If $x\neq y$, then $(x,y)$ and $(y,x)$ are two distinct elements in $\mathcal{H}_{N,m}$. However, they pair up to the same solution $x+y\mod m$ in $\mathcal{L}_{N,m}$. One can prove that, for a prime modulus $m$, the set $\mathcal{L}_{N,m}$ contains only about half of the possible residues modulo $m$. If $m$ is the product of many small primes, the cardinality of $\mathcal{L}_{N,m}$ is just a \emph{small fraction} of $m$. Additionally, the value of $p+q\pmod m$ is always in $\mathcal{L}_{N,m}$. As a consequence, these observations allow to deduce significant information about $S:=p+q$, the sum of the prime divisors of $N$. Here is a brief example: Let $N:=7909787=3823\cdot 2069$. For $r=5$, we have $N\equiv 2 \mod 5$ and $(1,2),(2,1),(3,4),(4,3)\in\mathcal{H}_{N,5}$. Consequently, $2,3\in\mathcal{L}_{N,5}$. Considering all primes up to $11$ and putting $m=2^2\cdot 3\cdot 5\cdot 7\cdot 11$, we observe that $\mathcal{L}_{N,m}$ contains only $40$ elements. Therefore, $S\pmod m$ is restricted to $40/4620=0.8\%$ of the residue classes modulo $m$.

In \cite{Hit1}, the sieve has been applied to improve the complexity of a babystep-giantstep routine modulo $N$. In the present paper, we will consider other possible applications of the idea. The remainder of the paper is structured as follows. \Cref{sec:prelim} concerns a preliminary result on Chinese remaindering. In \Cref{sec:sieve}, we will extend the hyperbolic sieve to general linear combinations $au+bv$ of co-divisors $u$ and $v$ of $N$. We will prove some useful properties and bounds for the cardinality of the related sieve sets. In \Cref{sec:fermat1}, we apply the sieve to Fermat's factorization method. While being one of the oldest and most basic techniques in the field, Fermat's algorithm is most efficient if it is applied to $N=uv$ with small divisor difference $\Delta:=|u-v|$. It is well known that the RSA cryptosystem is insecure if the difference of the prime factors of the public key $N$ is too small. The problem of small prime difference is frequently mentioned in standards, e.g. in ANSI \cite{Ans}, which requires a difference in the first 100 bits of the two prime factors of $N$. As a special purpose algorithm, Fermat's algorithm is thus still subject to research efforts and cost improvements. In \Cref{thm:main}, we refine the runtime complexity of Fermat's original procedure by a sub-exponential factor, while the space-complexity remains negligible. According to our knowledge of the literature, this yields the best known integer factorization bound in terms of the divisor difference $\Delta$. \Cref{sec:reduction} is then concerned with an extension of the hyperbolic sieve, leading to the reformulation of the integer factorization problem as an instance of the multiple choice subset-sum problem. As discussed in \Cref{rem:mcss}, this approach may have potential for further applications in the area of practical factorization. For instance, it is used in \Cref{sec:fermat2} to obtain an additional heuristic improvement of the bound in \Cref{thm:main} in return for higher space complexity. Finally, one of the more technical proofs of \Cref{sec:sieve} can be found in the appendix.

\section{Optimal Chinese remaindering \label{sec:prelim}}
In this short section, we recall an efficient way for computing sets via the Chinese Remainder Theorem that has first been presented in \cite{Hit1}. Let $m_1,\ldots,m_d\in\N$ be pairwise coprime and $\mathcal{A}_i\subseteq \Z_{m_i}$ for $i=1,\ldots,d$. Define $M:=\prod_{i=1}^d m_i$. As a result of the Chinese Remainder Theorem, we know that there is a unique set $\mathcal{A}\subseteq\Z_M$ induced by the sets $\mathcal{A}_i$. We may define this set as
\[
\mathcal{A}:=\left\{x\in\Z_M\text{ }\Bigg{|}\text{ } \exists (\alpha_1,\ldots,\alpha_d)\in\prod_{i=1}^d \mathcal{A}_i:x\equiv \alpha_i\mod m_i,\text{ } 1\leq i \leq d\right\}.
\]
We will now discuss an asymptotically optimal algorithm for determining the set $\mathcal{A}$. For $i=1,\ldots,d$, we set $\kappa_i:=|\mathcal{A}_i|$ and $\mathcal{A}_i=\{a_{i,1},\ldots,a_{i,\kappa_i}\}$. Furthermore, let $M_i:=Mc_i/m_i$, where $c_i:=(M/m_i)^{-1} \mod m_i$. Following the Gaussian approach to represent all elements of $\mathcal{A}$, we want to compute
\[
	a_{1,l_1}M_1+a_{2,l_2}M_2+\cdots+a_{d,l_d}M_{d} \mod M,
\]
where $l_i\in\{1,\ldots,\kappa_i\}$. The idea of the following procedure is to precompute a set of distances $\Delta$ for reaching all such elements, and then cycle through them in a clock-wise manner. It is worth to note that the space complexity of the mentioned precomputation step is
$
O(\sum_{i} \kappa_i),
$
which is about the same as the input size. This will be important for achieving negligible space complexity in our improvement of Fermat's factorization method in \Cref{sec:fermat1}. Let us now consider the algorithm.

\begin{algorithm}
	\emph{Input:} Pairwise coprime $m_i\in\N$, $\mathcal{A}_i\subseteq \Z_{m_i}$ and $\kappa_i$ for $i=1,...,d$, where $\kappa_1\geq \kappa_2\geq\ldots\geq \kappa_d$.
	
	\emph{Output:} The set $\mathcal{A}\subseteq\Z_M$ determined by the $\mathcal{A}_i$ due to the CRT.
	
	\begin{algorithmic}[1]
		\State Compute $M$ and $|\mathcal{A}|=\prod_{i=1}^d \kappa_i$. 
		\For {$i=1,\ldots,d$} \Comment{Precomputation  loop}
		\State Compute $M_i\mod M$ and set $r_i=0$. 
		\State $\Delta_{i,0}\gets (a_{i,1}-a_{i,\kappa_i})M_i\mod M$
		\For{$l_i=1,\ldots,\kappa_i-1$}
		\State $\Delta_{i,l_i}\gets (a_{i,l_i+1}-a_{i,l_i})M_i \mod M$ 
		\EndFor
		\EndFor
		\State $x_1\gets a_{1,1}M_1+...+a_{d,1}M_d \mod M$
		\For {$\nu=1,\ldots,|\mathcal{A}|-1$} \Comment{Main loop}
		\State Compute $r_1\gets r_1+1\mod{\kappa_1}$ and set $\mu=1$.
		\While {$r_\mu=0$}
		\State $\mu\gets \mu+1$
		\State $r_\mu\gets r_\mu+1\mod{\kappa_\mu}$
		\EndWhile
		\State $x_{\nu+1}\gets x_{\nu}+\sum_{i=1}^\mu\Delta_{i,r_i} \mod M$
		\EndFor
	\end{algorithmic}
	\label{alg:fastcrm}
\end{algorithm}

{\prop{\Cref{alg:fastcrm} is correct. Assuming $|\mathcal{A}|\geq \sum_{i=1}^d \kappa_i\cdot \log M$, its runtime complexity is bounded by $O(|\mathcal{A}|\cdot \log M)$ bit operations.}}

\begin{proof}
This is Theorem 4.2 in \cite{Hit1}.
\end{proof}

\section{Generalizing the hyperbolic sieve \label{sec:sieve}}
In the introduction, we have already discussed the sieve technique first established in \cite{Hit1}. We now give a more formal and general definition.

{\definition{Let $N,m,k\in\N$ such that $\gcd(Nk,m)=1$. We define the sets
\[
\mathcal{H}_{N,m}:=\{(x,y)\in\Z_m^2 :N\equiv xy \mod m\}
\]
and
\[
\mathcal{L}_{N,m,k}:=\{kx+y\mod m: (x,y)\in\mathcal{H}_{N,m}\}.
\]}\label{def:hyper}}

The set $\mathcal{H}_{N,m}$ is called a modular hyperbola and has the following properties:
\begin{enumerate}
\item{$|\mathcal{H}_{N,m}|=\phi(m)$.}
\item{If $N=uv$, then $(u\mod m, v \mod m)\in \mathcal{H}_{N,m}$.}
\end{enumerate}
Clearly, $N=uv$ implies $ku+v\pmod m\in\mathcal{L}_{N,m,k}$. 
In \cite{Hit1}, the set $\mathcal{L}_{N,m}:=\mathcal{L}_{N,m,1}$ has been used to search for the sum of the prime factors of semiprime numbers $N$. In the present paper, we are interested in general linear combinations of the factors of $N$. We now prove some useful properties of the corresponding sets.

{\lemma{Let $N,m,k\in\N$ such that $\gcd(Nk,m)=1$. For all $(a,b)\in\mathcal{H}_{k,m}$, we have
		\[
		\mathcal{L}_{N,m,k}=\{ax+by\mod m: (x,y)\in\mathcal{H}_{N,m}\}.
		\]}}

\begin{proof}
	Our assumptions imply that $\gcd(b,m)=1$. Hence, each pair $(x,y)\in\mathcal{H}_{N,m}$ corresponds uniquely to the pair $(\bar{x},\bar{y}):=(bx\pmod m,b^{-1}y\pmod m)\in\mathcal{H}_{N,m}$. As a result, the set on the right-hand side of the equation is equal to 
	\[
	\{a\bar{x}+b\bar{y}\mod m: (\bar{x},\bar{y})\in\mathcal{H}_{N,m}\}.
	\]
	It is easy to see that this proves the claim.
\end{proof}

As a result, the residues of all linear combinations $ap+bq$ of the divisors of $N$ with $ab\equiv k \mod m$ are in $\mathcal{L}_{N,m,k}$. In particular, this holds if $ab=k$.

{\lemma{Let $N,m,k\in\N$ such that $\gcd(Nk,m)=1$. If $s\in\mathcal{L}_{N,m,k}$, then $-s\in\mathcal{L}_{N,m,k}$.}\label{lem:props}}

\begin{proof}
	By definition, we have $s\equiv kx+y \mod m$ for some $(x,y)\in\mathcal{H}_{N,m}$. It follows that $N\equiv xy\equiv (-x)(-y) \mod m$, hence $-s\equiv k(-x)+(-y)\mod m$ is also in $\mathcal{L}_{N,m,k}$.
\end{proof}

{\lemma{
Let $N,k,e\in\N$ and $r$ be an odd prime such that $\gcd(Nk,r)=1$. Then
\[
|\mathcal{L}_{N,r,k}|=(r+(k^{-1}N|r))/2,
\] 
where $(\cdot|\cdot)$ is Legendre's symbol.}\label{lem:card}}

\begin{proof}
Our approach is to identify those pairs in $\mathcal{H}_{N,r}$ which yield the same element in $\mathcal{L}_{N,r,k}$. Let $(x_0,y_0)$ and $(x_1,y_1)$ be distinct elements of $\mathcal{H}_{N,r}$ such that they pair up to one element in $\mathcal{L}_{N,r,k}$, namely by satisfying $kx_0+y_0\equiv kx_1+y_1 \mod r$. One easily observes that this is equivalent to
\[
kx_0^2-x_0(kx_1+y_1)+x_1y_1\equiv 0 \mod r.
\]
Considering this as quadratic congruence in the variable $x_0$, it follows that $x_0=x_1$ or $x_0=k^{-1}y_1\pmod r$. The first case implies that $y_0=y_1$ and hence is not possible due to our assumption. We derive that the elements $(x_1,y_1)$ and $(k^{-1}y_1 ,kx_1)$ of $\mathcal{H}_{N,r}$ pair up to one element in $\mathcal{L}_{N,r,k}$, and that these are the only ones that do so. Still, two such pairs can be equal. This happens if and only if $y_1\equiv kx_1\mod r$, which yields
\[
N\equiv x_1y_1\equiv kx_1^2 \mod r.
\]
If $k^{-1}N$ is a quadratic nonresidue modulo $r$, the latter case is not possible and we derive that there are exactly $(r-1)/2$ elements in $\mathcal{L}_{N,r,k}$. However, if $(k^{-1}N|r)=1$, there are two distinct solutions to $x^2\equiv k^{-1}N\mod r$. We conclude that there are exactly $(r-3)/2+2=(r+1)/2$ elements in $\mathcal{L}_{N,r,k}$.	
\end{proof}

The following statement describes the behavior of the cardinality of $\mathcal{L}_{N,m,k}$ for $m$ being an odd prime power. The proof is rather involved, technical and has no strong relevance to the main contents of the paper. It can be found in the appendix.
	
{\theorem{Let $N,k,e\in\N$ and $r$ be an odd prime such that $\gcd(Nk,r)=1$. Then
for $e\geq 1$, we have
\[
|\mathcal{L}_{N,r^{e+1},k}|=
 \begin{cases}
 |\mathcal{L}_{N,r^e,k}|\cdot r     & \mbox{if } (k^{-1}N|r)=-1,\\
|\mathcal{L}_{N,r^e,k}|\cdot r-2^{e\pmod 2}\cdot (r-1)     & \mbox{if } (k^{-1}N|r)=1.
   \end{cases}
\]}\label{thm:sievepower}}

Assume that we know the prime factorization of $m$. The following lemma proves that we may compute the sets $\mathcal{L}_{N,m,k}$ by applying \Cref{alg:fastcrm}.

{\lemma{Let $N,k\in\N$ and $m_1,\ldots,m_d\in\N$ be pairwise coprime. Furthermore, define $m:=m_1\cdots m_d$. If $\gcd(Nk,m)=1$, then $\mathcal{L}_{N,m,k}\subseteq \Z_m$ is the unique set determined by the $\mathcal{L}_{N,m_i,k}$ due to the Chinese Remainder Theorem (CRT).}\label{lem:crt}}

\begin{proof}
Consider the map
\begin{align*}
&\psi:\mathcal{L}_{N,m,k}\rightarrow \prod_{i=1}^{d} \mathcal{L}_{N,m_i,k},\\
x\mapsto (x&\mod m_1,\ldots,x\mod m_d).
\end{align*} 
We first show that $\psi$ is well-defined. Let $x\in\mathcal{L}_{N,m,k}$ and $i\in\{1,\ldots,d\}$ be arbitrarily chosen. By definition, there exists at least one pair $(x',y')\in\mathcal{H}_{N,m}$ such that $x\equiv kx'+y' \mod m$. Clearly, $(x'\mod m_i, y'\mod m_i)$ is an element of $\mathcal{H}_{N,m_i}$. We deduce that the element $x \mod m_i=kx'+y' \mod m_i$ is in $\mathcal{L}_{N,m_i,k}$, which we wanted to prove.

As a result of the CRT, $\psi$ is injective. We prove that this map is also surjective. Let $(s_1,\ldots,s_d)\in \prod_{i=1}^{d} \mathcal{L}_{N,m_i,k}$ be arbitrary. For every $i\in\{1,\ldots,d\}$, there exists a pair $(u_i,v_i)$ in $\mathcal{H}_{N,m_i}$ such that $s_i\equiv ku_i+v_i\mod m_i$. By the CRT, we exhibit $u$ and $v$ such that
\[
u\equiv u_i \mod m_i \text{ and } v\equiv v_i \mod m_i
\]
and, hence, $N\equiv uv \mod m_i$ for every $i$. This implies that $N\equiv uv \mod m$, and $s:=ku+v\mod m\in\mathcal{L}_{N,m}$. It is easy to see that $\psi(s)=(s_1,\ldots,s_d)$. Therefore, $\psi$ is bijective and the statement follows. 
\end{proof}

Finally, we give a bound for the cardinality of the sieve set in the case where $m$ is the product of several small prime factors.
{\theorem{
Let $N,k\in\N$, $B\in\R^+$ and set $m=\prod_{2<r\leq B} r$ for primes $r$. Assume that $\gcd(Nk,m)=1$. Then
\[
|\mathcal{L}_{N,m,k}|\in O\left(\frac{m\log B}{2^{\omega(m)}}\right),
\]
where $\omega(m)$ is the number of prime divisors of $m$.
}\label{thm:sievebound}}

\begin{proof}
We first note that 
\[
\left(\prod_{r\leq B}\frac{r+1}{r}\right)\in\Theta(\log B),
\] 
which is an immediate consequence of Mertens' theorem (\cite[Thm. 5.13]{sho}) and the fact that
\[
\left(\prod_{r\leq B}\frac{r+1}{r}\right)\left(\prod_{r\leq B}\frac{r-1}{r}\right)=\prod_{r\leq B} \left(1-\frac{1}{r^2}\right)\in\Theta(1).
\]
As a result of \Cref{lem:crt}, we derive $|\mathcal{L}_{N,m,k}|=\prod_{2<r\leq B}|\mathcal{L}_{N,r,k}|$. Note that \Cref{lem:card} implies $|\mathcal{L}_{N,r,k}|\leq (r+1)/2$ for odd primes $r$. We conclude that
\[
\frac{|\mathcal{L}_{N,m,k}|}{m}\leq \frac{4}{3}\cdot\prod_{r\leq B}\frac{(r+1)/2}{r}\in O\left(\frac{\log B}{2^{\omega(m)}}\right).
\]
\end{proof}

{\remark{Let $f\in\Z[X,Y]$. A further, natural generalization of the hyperbolic sieve would be to consider sets of the form
		\[
		\mathcal{L}_{N,m,f}:=\{f(x,y)\pmod m: (x,y)\in\mathcal{H}_{N,m}\}.
		\]
		Clearly, we have $f(p,q)\pmod{m}\in\mathcal{L}_{N,m,f}$. Also, if $f$ is symmetric, then the fundamental theorem of symmetric polynomials implies that $|\mathcal{L}_{N,m,f}|\leq |\mathcal{L}_{N,m}|$. However, the increased average size of values of $f(p,q)$ for non-linear polynomials $f$ makes it more difficult to make proper use of the information provided by the sieve. On the other hand, $\mathcal{L}_{N,m,f}$ may be significantly smaller than $\mathcal{L}_{N,m}$ for polynomials $f$ of higher degree. For example, taking $f=(X+Y)^n$ and primes $r$ with $n\mid r-1$, we have $|\mathcal{L}_{N,r,f}|\leq (r-1)/n+1$. This fact might also be interesting in light of the approach discussed in \Cref{sec:reduction}.
}}
	
\section{Fermat factorization: Sieve improvement}\label{sec:fermat1}
Fermat's factorization algorithm has first been described by Pierre de Fermat in 1643. Its main 
idea is still used in modern factorization algorithms and is based on the representation of $N$ as 
the difference of two squares. Such representation allows to write
$
N=X^2-Y^2=(X+Y)(X-Y),
$
hence $X+Y$ often yields a proper factor of $N$. In a generalization of Fermat's original procedure from 1895, Lawrence \cite{Law} applied the very same idea to multiples of the number $N$. Let us discuss this approach in greater detail. In order to find suitable square numbers, we may use the fact that $(au+bv)^2-(au-bv)^2=4abN$ holds for all $a,b,u,v\in\N$ with $N=uv$. Note that this equality also implies
\[
(au-bv)^2=(au+bv)^2-4abN=(au+bv-2\sqrt{abN})(au+bv+2\sqrt{abN}),
\]
hence
\[
0\leq au+bv-2\sqrt{abN}=\frac{(au-bv)^2}{au+bv+2\sqrt{abN}}<\frac{(au-bv)^2}{4\sqrt{abN}}.
\]
Let $L:=\lceil 2\sqrt{abN}\rceil$. We want to find $z\in\N$ such that $L+z=au+bv$. Considering the bound shown above, we know that 
$z<(au-bv)^2/(4\sqrt{abN})$. Since $N$ may have several nontrivial decompositions $N=uv$, we are interested in the minimal upper bound for $z$ with regards to $a$ and $b$ that will allow us to find a nontrivial divisor.
{\definition{Let $a,b\in\N$ and $N$ be a composite number. We define
\begin{equation}
\Lambda_{a,b}:=\min\left\{\frac{(au-bv)^2}{4\sqrt{abN}}\,\Big{|}\,(u,v)\in\N^2\,:\, N=uv \wedge u,v\neq 1\right\}.
\end{equation}
}\label{def:zbound}}
For each value $i<\Lambda_{a,b}$, we have to check if $(L+i)^2-4abN$ is a square number. Fermat's original algorithm for $a=b=1$ checks all possible values $i<\Lambda_{1,1}$, hence runs in time $\widetilde{O}(\Lambda_{1,1})$. This runtime complexity is small if there exist co-divisors $u,v$ of $N$ such that $|u-v|$ is small. We may thus consider Fermat's approach as special purpose factorization algorithm for factoring numbers with small divisor difference. As such, the algorithm has been improved several times, e.g. in \cite{McK}, \cite{SK} and \cite{WTS}. Most of these improvements concern the use of a modular sieve for reducing the number of candidates $i$ for $z$. If $m$ is some small natural number, we only have to check those $i<\Lambda_{1,1}$ for which $(i+\lceil 2\sqrt{N}\rceil)^2-4N$ is a square modulo $m$. For example, if $N=7909787$ and $m=20$, the value of $z\pmod m$ has to be in the set $\{2,3,7,8,12,13,17,18\}$. As a result, this sieve decreases the overall number of values $i$ we need to consider. However, we have to precompute and store a table of the residues modulo $m$, distinguishing between squares and non-squares. Some authors propose the use of certain moduli, like $m=11,63,64$ and $65$. The number of candidates for $z$ surviving this sieve reduces the runtime complexity by a large, but only numerical factor (see \cite[p.149]{Rie}, \cite[p.125]{Wag}).

One of the main contributions of this paper is to provide a much more substantial reduction of the runtime complexity of Fermat's factorization algorithm and its generalization by Lawrence. The space complexity is still negligible and implementation is quite straightforward. Instead of the sieve described in the last paragraph, we will apply the hyperbolic sieve to obtain an \emph{asymptotic} runtime improvement by a \emph{sub-exponential} factor. We start by observing that the hyperbolic sieve is at least as strong as the sieve described above, as the following result shows.

{\lemma{Let $N,m,k\in\N$ such that $\gcd(Nk,m)=1$. Then
		\[
	\mathcal{L}_{N,m,k}\subseteq \{x\in\Z_m:x^2-4kN \text{ is a square modulo $m$}\}.
	\]
If $m$ is odd, the sets are equal.		
}\label{lem:qres}}

\begin{proof}
	Let $x=ka+b\in\mathcal{L}_{N,m,k}$ be arbitrary, where $(a,b)\in\mathcal{H}_{N,m}$. We have
	\[
	x^2-4kN\equiv (ka+b)^2-4kab\equiv (ka-b)^2 \mod m.
	\]
	So it follows that $x^2-4kN$ is a square modulo $m$,
	which proves the first claim. We now show that, for odd $m$, we have equality. Let $x\in\Z_m$ such that $x^2-4kN$  is a square modulo $m$. It remains to prove that $x\in\mathcal{L}_{N,m,k}$. There exists $y\in\Z$ such that 
	\[
	x^2-y^2\equiv (x+y)(x-y)\equiv 4kN \mod m.
	\]
	Hence, there is $d\in \Z_m^*$ with $x+y\equiv 2d \mod m$ and $x-y\equiv (4kN)(2d)^{-1}\mod m$, implying that
	\[
	x\equiv k(Nd^{-1})+d\mod m.
	\]
	So we have shown that $x$ is of the form $ka+b\pmod m$ with $ab\equiv N\mod m$, hence $x\in\mathcal{L}_{N,m,k}$.
\end{proof}

Let $a,b,m\in\N$ such that $\gcd(abN,m)=1$. In addition, let $k:=ab \pmod m$ and $\Lambda_{a,b}$ like in \Cref{def:zbound}. From the results in \Cref{sec:sieve}, it follows that the value $z$ we are looking for satisfies
$
z\equiv s-L \mod m
$
for $L:=\lceil 2\sqrt{abN}\rceil$ and some $s\in\mathcal{L}_{N,m,k}$. We hence define
\[
\bar{\mathcal{L}}_m:=\{s-L\pmod m: s\in\mathcal{L}_{N,m,k}\}.
\]
For every $i\in\bar{\mathcal{L}}_m$, we then need to check if $(i+L)^2-4abN$ is a square. If we have $(i_0+L)^2-4abN=y^2$, we compute $\gcd(i_0+L-y,N)$. Note that if $m$ is larger than $\Lambda_{a,b}$, then $z\in\bar{\mathcal{L}}_m$ and we obtain a nontrivial divisor of $N$ as a result. 

In order to make our improvement as efficient as possible, we have to deal with two questions: How should we choose the modulus $m$? And what is the best way to compute the set $\bar{\mathcal{L}}_m$?

According to \Cref{thm:sievebound}, we should take $m$ as the product of small primes up to a certain bound $B$. Ideally, we would choose $m=\prod_{2<r\leq B} r$ with $B$ minimal such that $m>\Lambda_{a,b}$. Since we do not know $\Lambda_{a,b}$, we will need to pick a bound $\Lambda$ to run our algorithm with. This bound determines the search interval for $z$ and, hence, the overall runtime complexity. In any case, we will see that working with such highly composite moduli $m$ allows for a substantial reduction of the runtime. In addition, \Cref{alg:fastcrm} allows to compute the elements of the set $\bar{\mathcal{L}}_m$ \emph{consecutively} and in \emph{asymptotically optimal time}. This way, we do not have to store tables of possible candidates and may keep the space complexity of the algorithm as small as possible.

\begin{algorithm}
	\emph{Input:} A composite, odd number $N$ and a bound $\Lambda\in\N$
	
	\emph{Output:} A proper factor of $N$ or `$\Lambda$ is too small'
	
	\begin{algorithmic}[1]
		\State Let $m=\prod_{2<r\leq B} r$ with $B$ minimal such that $m>\Lambda$.
		\State Compute $\gcd(N,m)$. If a factor of $N$ is found, print and stop.
		\State Choose $a,b\in\N$, $a,b<N$ with $\gcd(ab,mN)=1$. 
		\State Let $L:=\lceil 2\sqrt{abN}\rceil$ and $k:= ab\mod m$.
		\For {$r\mid m$} 
		\State Compute $\bar{\mathcal{L}}_{r}:=\{s-L\pmod m: s\in\mathcal{L}_{N,r,k}\}$. 
		\EndFor
		\State Run \Cref{alg:fastcrm} to compute $\bar{\mathcal{L}}_m$ from the sets $\bar{\mathcal{L}}_{r}$.
		\For {$i\in\bar{\mathcal{L}}_m$}
		\If{$(i+L)^2-4abN = y^2$ for some $y\in\N$} print $\gcd(i+L-y,N)$.
		\EndIf
		\EndFor
		\State If no proper factor is found, return `$\Lambda$ is too small'.
	\end{algorithmic}
	\label{alg:sieveimp}
\end{algorithm}

{\prop{		
If 
$
\Lambda>\Lambda_{a,b},
$
then \Cref{alg:sieveimp} returns a proper factor of $N$. Ignoring terms that are polynomial in $\log N$ and $\log \Lambda$, the runtime complexity is given by 
\[
\widetilde{O}(\Lambda\cdot \exp(-C\log\Lambda/\log\log\Lambda))
\]
bit operations, where $C=(1+o(1))\log 2$. 
}\label{prop:sieveimp}}

\begin{proof}
	We first prove correctness. Assume that $\Lambda>\Lambda_{a,b}$. In the Steps 1--4, we choose our parameters such that $\gcd(abN,m)=1$ and $\gcd(ab,N)=1$ holds. In \Cref{lem:crt}, we have seen that the set $\mathcal{L}_{N,m,k}$ is induced by the sets $\mathcal{L}_{N,r,k}$ due to the Chinese Remainder Theorem. One easily observes that the same is true for the sets $\bar{\mathcal{L}}_m$ and $\bar{\mathcal{L}}_{r}$. Hence, we may compute the elements of $\bar{\mathcal{L}}_m$ by applying \Cref{alg:fastcrm} in Step 7. It remains to show that we will obtain a proper factor at some point in the final for-loop in Step 9. From $\Lambda>\Lambda_{a,b}$, it follows that there exists a nontrivial composition $N=uv$ of $N$ such that
	\[
	\Lambda>\frac{(au-bv)^2}{4\sqrt{abN}}>z,
	\]
	where $z+L=au+bv$ (see the discussion at the beginning of \Cref{sec:fermat1}). In addition, $z=z\pmod m$ is an element of the set $\bar{\mathcal{L}}_m$. Hence, for the run $i=z$ of the for-loop, we will find $y=au-bv$ and compute the GCD of $i+L-y=2bv$ and $N$. Since $N$ is odd and $\gcd(ab,N)=1$, this yields a proper factor of $N$.
	
	Let us now consider the runtime complexity of the algorithm. We first discuss the magnitude of $B$ in Step 1. For Chebyshev's theta function $\vartheta(x):=\sum_{r\leq B} \log r$, it follows that
	\[
	e^{\vartheta(B)}=\prod_{r\leq B} r =2m\in \Theta(\Lambda).
	\]
	From the Prime Number Theorem, it follows that $\vartheta(x)-x\in o(x)$ (see \cite[p.80]{RosSch}). We conclude that $B=(1+o(1))\log \Lambda$ for increasing values of $\Lambda$. It immediately follows that the precomputations in the Steps 1--5 are polynomial in $\log N$ and $\log \Lambda$. For determining the $\bar{\mathcal{L}_r}$ in Step 6, we just use \Cref{def:hyper}.  Note that the cardinalities of these sets are bounded by $B\in O(\log \Lambda)$. Let us now consider the cost for the Steps 7--10. Since the computations in Step 9 are all polynomial in $\log N$, the main cost depends on the size of $\bar{\mathcal{L}}_m$. For the prime-counting function $\pi$, \Cref{thm:sievebound} yields
	\[
	\widetilde{O}\left(\frac{m}{2^{\omega(m)}}\right)=\widetilde{O}\left(\frac{\Lambda}{2^{\pi(B)}}\right).
	\]
The Prime Number Theorem states that $\pi(B)\approx B/\log B$. Indeed, we have the explicit lower bound of $B/\log B<\pi(B)$ for $B\geq 17$ (see \cite[p.69]{RosSch}). Hence, for $C=(1+o(1))\log 2$, we obtain
	\begin{align*}
		O(2^{-\pi(B)})= O\left(2^{-\frac{B}{\log B}}\right)&=O\left(2^{-\frac{(1+o(1))\log\Lambda}{\log\left(\log\Lambda\right)}}\right)\\
		&= O\left(\exp(-C\log\Lambda/\log\log\Lambda\right),
	\end{align*}
	which proves the claim. 
\end{proof}

{\remark{The rigorously proved runtime complexity represents the worst case. The algorithm may find a factor faster if we are lucky and $z$ is actually among the earlier elements tested in the for-loop in Step 9. In addition, we want to stress that there is no need to compute and store the complete set $\bar{\mathcal{L}}_{m}$. We may use \Cref{alg:fastcrm} to compute its elements consecutively and run Step 9 of \Cref{alg:sieveimp} for each of them.  If no factor is found, we delete the element and go to the next one. In the precomputation loop of \Cref{alg:fastcrm}, we only need the sets $\bar{\mathcal{L}}_{r}$. Since their cardinalities are bounded by $B\in O(\log \Lambda)$, the required space is negligible.
}\label{rem:space}}
\vspace{12pt}

Unfortunately, we cannot apply \Cref{alg:sieveimp} to values for $a$ and $b$ that are divisible by primes smaller than $B$, at least not without increased runtime costs. However, we can always apply it to the case of Fermat's original factoring method, namely to $a=b=1$. In this regard, we now prove a deterministic and rigorous factorization bound in terms of the minimal divisor difference.

{\theorem{Let $N$ be a composite number and $\Delta:=\Lambda_{1,1}$ as in \Cref{def:zbound}. Ignoring terms that are polynomial in $\log N$ and $\log \Delta$, there exists a deterministic algorithm that computes a proper factor of $N$ in 
\[
\widetilde{O}(\Delta\cdot \exp(-C\log\Delta/\log\log\Delta))
\]
bit operations, where $C=(1+o(1))\log 2$. The space complexity is negligible.
}\label{thm:main}}

\begin{proof}
We run \Cref{alg:sieveimp}  with $a=b=1$ and for increasing values of the modulus $m$, starting with $m=2,2\cdot 3, 2\cdot 3\cdot 5,\ldots$, until we find a proper factor. We succeed as soon as $m>\Delta$. Hence, a (rough) estimate for the total cost is given by the runtime of \Cref{prop:sieveimp} multiplied with a factor of $\log \Delta$. However, this factor is negligible. The space complexity is discussed in \Cref{rem:space}.
\end{proof}

In \Cref{alg:sieveimp}, the modulus $m$ is a product of small odd primes. Based on \Cref{thm:sievepower}, we will now discuss some practical improvements by constant factors. Depending on $N$ and $k$, we may want to include powers of $2$ and other small prime powers in order to keep the ratio between $|\bar{\mathcal{L}}_m|$ and $m$ as small as possible. We illustrate this with the following example.

{\example{Let $N=17344343992304993085649094809$. We want to search for a linear combination of factors of $N$ with $a=b=1$. We note that $N\equiv 1 \mod 8$ and $(N|r)=1$ for $r=3,5,7$. In order to build $m$, we consider the following:
		\begin{itemize}
			\item{For $r=2$, our computations (see \Cref{rem:twopower} in the appendix) show that
				$|\bar{\mathcal{L}}_{m}|=2,4,4,6,8,14,24$ for $m=2^4,2^5,2^6,2^7,2^8,2^9,2^{10}$.
			\item{For $r=3$, \Cref{thm:sievepower} implies 
				$|\bar{\mathcal{L}}_{m}|=2,2,4,8,22$ for $m=3,3^2,3^3,3^4,3^5$, respectively.}
			\item{For $r=5$, \Cref{thm:sievepower} implies 
				$|\bar{\mathcal{L}}_{m}|=3,7,31$ for $m=5,5^2,5^3$.}
			\item{For $r=7$, \Cref{thm:sievepower} implies 
				$|\bar{\mathcal{L}}_{m}|=4,16,106$ for $m=7,7^2,7^3$.}
			}
		\end{itemize}
	Our goal is to keep $|\bar{\mathcal{L}}_m|/m$ as small as possible. Adding higher powers of the primes $2$, $3$ and $5$ appears to benefit this proportion. In our example, going from $3$ to $3^2$ in $m$ increases $m$ by a factor of $3$ while leaving $|\bar{\mathcal{L}}_m|$ the same. Adding $5^2$ to $m$ is also beneficial. Starting with $7$, however, we always benefit more from adding the next highest prime than from adding another power of the same prime. As we can see above, going from $7$ to $7^2$ in $m$ adds a factor of $(16/4)/7=4/7$ to $|\bar{\mathcal{L}}_m|/m$, while including $11$ in $m$ adds a factor of either $5/11$ or $6/11$ (see \Cref{lem:card}), both being smaller than $4/7$.
	
	Of course, there are also other considerations that may affect how to build the best modulus $m$. For instance, assume that we are already close to our target bound $\Lambda$. Even if adding a larger prime might be better for the proportion, it may make much more sense to just add another small factor like $2$ or $3$. For this reason, we will not give an explicit algorithm for building $m$, as the best choice may depend on specific circumstances.		
		
		For factoring $N$, we take $m=2^8\cdot 3^3\cdot 5^2\cdot 7\cdot 11\cdot 13\cdot 17\cdot 19=55870214400$. Then $\bar{\mathcal{L}}_m$ has $1935360$ elements, and \Cref{alg:sieveimp} returns
		\[
		N=pq=129411310904131\cdot 134024946282739.
		\]
Note that $p+q-\lceil 2\sqrt{N}\rceil=40403063803$. So our algorithm only needs to check at most $0.00479\%\approx 1935360/40403063803$ of the candidates that would have been checked in Fermat's original procedure.
}\label{exa:powers}}

\section{Reducing integer factorization to multiple choice subset-sum}\label{sec:reduction}
We now describe a new research direction based on the hyperbolic sieve. Let $N$ be the number we want to factorize and $m\in\N$ coprime to $N$. For the sake of simplicity, we assume that $N=pq$ is a semiprime number and consider the sieve set $\mathcal{L}_{N,m}=\mathcal{L}_{N,m,1}$ for finding the sum $S=p+q$ of the prime factors of $N$. In the previous subsection, we chose $m$ as a product of primes in accordance with \Cref{thm:sievebound}. In addition, we tried to take $m$ as small as possible, just quite large enough to retrieve enough information about $S$ (see the proof of \Cref{thm:main}). The reason is simple: the smaller $m$ is, the less elements are in $\mathcal{L}_{N,m}$ that have to be checked as candidates for $S$. In \Cref{exa:powers}, the inclusion of all primes up to $19$ was sufficient to find $S$. While we may have easily computed $\mathcal{L}_{N,23}$ and larger sieve sets, we did not need to consider them. On the other hand, these larger sieve sets contain information about $S$ that is (almost) freely available to us, and it seems like there is room for further improvement by finding a way to make use of them. In this section, we will discuss such an approach in which we choose $m$ large. Then $S$ will be one of the smallest elements in $\mathcal{L}_{N,m}$.

The idea is based on the fact that $S$ is the only number which is an element of $\mathcal{L}_{N,m}$ for \emph{every} $m\geq S+1$. In other words, it holds that
$
\{S\}=\bigcap_{m\geq S+1}\mathcal{L}_{N,m}.
$
Assume that $p$ and $q$ are balanced and hence satisfy $p+q\leq 3\sqrt{N}$, which is common for the prime factors of semiprime RSA moduli. For $\bar{s}:=S-\lceil 2\sqrt{N}\rceil$, it easily follows that 
$
0\leq \bar{s} \leq \sqrt{N}.
$
Let $m$ be a product of the primes up to a small bound $B$ just like in \Cref{thm:sievebound}. Similar to \Cref{alg:sieveimp}, let $B$ be minimal such that $m>\Lambda$, where $\Lambda:=\lceil \sqrt{N}\rceil$. We already know that, among other elements, $\bar{s}$ is in 
$
\bar{\mathcal{L}}_m:=\{(x-\lceil 2\sqrt{N}\rceil) \pmod m: x\in\mathcal{L}_{N,m}\}.
$
We now continue to multiply $m$ with the first few prime numbers exceeding $B$, and denote the resulting modulus by $M$. By computing the corresponding sets $\bar{\mathcal{L}}_{M}$, one observes that the cardinality of  
\[
\bar{\mathcal{L}}_m \cap \bar{\mathcal{L}}_M =\{s\in\bar{\mathcal{L}}_M:s\leq \sqrt{N}\}
\]
decreases drastically with every additional prime added to $M$. This behavior may be explained by our previous results. Let $r$ be prime such that $r\nmid m$ and let $s\leq\sqrt{N}<m$ be an element of $\bar{\mathcal{L}}_m$. Lemma \ref{lem:crt} yields that $s$ is an element of $\bar{\mathcal{L}}_{mr}$ if and only if $s\pmod r\in\bar{\mathcal{L}}_r$. However, $\bar{\mathcal{L}}_r$ contains only about half of the possible residues modulo $r$. Therefore, one would expect that the number of elements $s$ in  $\bar{\mathcal{L}}_{M}$ satisfying $s\leq\sqrt{N}$ is \emph{cut about half} with every additional prime, until only $\bar{s}$ is left. So considering \Cref{thm:sievebound}, we would expect 
\begin{align}\label{eq:extbound}
|\{s\in\bar{\mathcal{L}}_M:s\leq \sqrt{N}\}|\in O\left(\frac{|\mathcal{L}_m|}{2^{\omega(M/m)}}\right)\subseteq O\left(\frac{\Lambda \log B}{2^{\omega(M)}}\right).
\end{align}
Depending on the desired cardinality of this set, it makes sense to consider different sizes for $M$. If we take $\omega(M)$ as large as $\log N$, we may certainly assume to have isolated $\bar{s}$ as the minimum of $\bar{\mathcal{L}}_M$. However, since we have not provided a rigorous proof for this statement, (\ref{eq:extbound}) will be used as heuristic assumption in the remainder of this paper.

\vspace{12pt}
For the sake of clarity, we will now suppose that $\bar{s}$ is the only element of $\bar{\mathcal{L}}_{M}$ satisfying $\bar{s}\leq\sqrt{N}$ and that the number $k:=\omega(M)$ of distinct prime factors of $M$ is in $O(\log N)$. Let $M=\rho_1\cdots \rho_k$, where the $\rho_i$ are either primes or powers of primes. As a consequence of \Cref{lem:crt}, we know that $\bar{\mathcal{L}}_{M}$ is determined by the $\bar{\mathcal{L}}_{\rho_i}$ due to the CRT. Just like in \Cref{sec:prelim}, we use the Gaussian approach to represent all elements in $\bar{\mathcal{L}}_{M}$. For $\gamma_i:=(M/\rho_i)^{-1}\pmod {\rho_i}$ and $M_i:=M\gamma_i/\rho_i$, every element in $\bar{\mathcal{L}}_{M}$ can be written as
\begin{align}\label{eq:crem}
	a_{1,j_1}M_1+a_{2,j_2}M_2+\ldots+a_{k,j_k}M_k\mod M,
\end{align}
where $a_{i,j_i}\in\bar{\mathcal{L}}_{\rho_i}$ with $j_i\in\{1,\ldots,\kappa_i\}$ for $|\bar{\mathcal{L}}_{\rho_i}|=\kappa_i$. Our ultimate goal is to find the choice of $(a_{1,j_1},\ldots,a_{k,j_k})$ for which the value of (\ref{eq:crem}) equals $\bar{s}$. According to our assumptions, however, this coincides with the task to find those $a_{i,j_i}$ for which (\ref{eq:crem}) is minimal. So in order to find $\bar{s}$ and factorize $N$, all we have to do is to minimize the expression (\ref{eq:crem}). 
This problem is strongly related to the \emph{multiple choice subset-sum problem} (MCSS), which is discussed in \cite[Sec. 11.10.1]{KelPfePis}. Before we explain this relation, let us first give a general description of the MCSS problem.
We consider the classes $N_i$ for $i=1,\ldots,k$, where each class contains the weights $w_{i,1},\ldots,w_{i,\kappa_i}$. The problem is to select exactly one weight from each class such that the total weight sum is maximized without exceeding a capacity $c$:
\begin{align}\label{eq:mcss}
	\text{maximize } &\sum_{i=1}^k\sum_{j=1}^{\kappa_i}w_{i,j}x_{i,j}\\\nonumber
	\text{subject to } &\sum_{i=1}^k\sum_{j=1}^{\kappa_i}w_{i,j}x_{i,j}\leq c,\\\nonumber
	& \sum_{j=1}^{\kappa_i}x_{i,j} = 1\text{, where } i=1,\ldots,k,\\\nonumber
	& x_{i,j}\in\{0,1\}\text{, where } i =1,\ldots,k\text{ and } j\in\{1,\ldots,\kappa_i\}.
\end{align}
Here, $x_{i,j}=1$ if item $j$ was chosen in class $N_i$. Other formulations of the problem concern minimization instead of maximization, approximation from both sides or even the exact representation of a specific target value as such sum. Recently, modular versions of  the problem have been investigated in \cite{AxiBacJin}. 

We are now ready to show how to reduce integer factorization to a MCSS instance. We will actually show two reductions, the first to the maximization problem formulated in (\ref{eq:mcss}), and the second to an exact target sum MCSS.

\textbf{Reduction to maximization MCSS (\ref{eq:mcss})}. We show that the structure of our minimization problem in (\ref{eq:crem}) may be translated into a MCSS problem instance. For $i=1,\ldots,k$, define
$
N_i:=\{a_{i,j_i}M_i \mod M: j_i=1,\ldots,\kappa_i\}.
$
We immediately see that all weights in $N_i$ are divisible by $M/\rho_i$. Note that the corresponding sums (\ref{eq:mcss}) never exceed $kM$. Given these classes, the modular version of the problem may be split into $k$ MCSSs as formulated above, namely one for each choice of capacity $c=\sqrt{N},M+\sqrt{N},\ldots, (k-1)M+\sqrt{N}$. Alternatively, we may also define the class $N_{0}:=\{0,M,2M,\ldots (k-1)M\}$, leading to only one MCSS with $c=(k-1)M+\sqrt{N}$. 

\textbf{Reduction to an exact target MCSS}. In addition, we show that we can also reduce to a MCSS problem instance with an exact target sum. We let $M=U\cdot V$ and make sure that $U,V >\sqrt{N}\geq \bar{s}$, where $U=\rho_1\cdots \rho_k$ and $V=\rho_{k+1}\cdots \rho_l$. Using the notation of (\ref{eq:crem}), we then have
\[
a_{1,j_1}U_1+\ldots+a_{k,j_k}U_k-j_0U\equiv 
b_{k+1,r_{k+1}}V_{k+1}+\ldots+b_{l,r_l}V_l \pmod V,
\]
where the $a_{i,j_i}$ and $b_{i,r_i}$ correspond to the residues of $\bar{s}$ in the sets $\bar{\mathcal{L}}_{\rho_i}$. The value $j_0\in\{0,\ldots,k\}$ is chosen  in a way such that the term on the left-hand side of the congruence is actually equal to $\bar{s}$. We then define the classes $N_0:=\{0,U,\ldots,kU\}$ and 
\[
N_i:=\{aU_i \mod V: a\in\bar{\mathcal{L}}_{\rho_i}\}
\]
for $i=1,\ldots,k$, and 
\[
N_i:=\{-bV_i \mod V: b \in\bar{\mathcal{L}}_{\rho_i}\}
\]
for $i=k+1,\ldots,l$. We hence obtain a MCSS instance with the exact target sum $0$ and 
significantly smaller weights, namely the size of $V$. \Cref{lem:crt} implies that $\bar{s}$ 
corresponds to a solution of the related MCSS instance, and if $U$ and $V$ are large enough, this 
is the only solution this problem will have. 
\vspace{12pt}

{\remark{		
		While general knapsack and subset-sum problems are known to be NP-complete and therefore unlikely to be solved in polynomial-time, a variety of algorithms has been developed for application, some of which with great success in specific settings. The employed methods range from branch-and-bound algorithms \cite{SinPosKol}, dynamic programming solutions and polynomial-time approximation schemes \cite[Sec. 4.1 and 4.5]{KelPfePis} to lattice reduction \cite{SchEuc} and even machine learning \cite{GuCui}. For an overview on methods developed specifically for the multiple-choice variant of subset-sum, see \cite[Sec. 11]{KelPfePis}.
		
		As mentioned, the (multiple choice) subset-sum problem is NP-complete and integer 
		factorization is in NP. While the existence of a reduction from the latter to the first 
		problem is thus not a new discovery, the goal of this section was to give an explicit and 
		tangible reformulation of integer factorization as a specific additive minimization 
		problem. The resulting reduction and its possible variants have the potential to introduce 
		new approaches to integer factorization. In particular, we plan on exploring the 
		application of lattice reduction techniques like LLL and BKZ. For now, we demonstrate the 
		promising nature of the reduction by showing that it allows for further improvement of 
		Fermat's factorization method.

}\label{rem:mcss}}

\section{Fermat factorization: Time-space tradeoff}\label{sec:fermat2}
Based on the reduction considered in the previous section, we now discuss a second factorization algorithm. In some sense, it is a time-space tradeoff for our improvement of Fermat's factorization algorithm (see \Cref{sec:fermat1}). With the notation of \Cref{sec:reduction} in mind, we consider only two classes in our MCSS and, hence, set $k=2$. Let $0<B_1<B_2$, $m_1=\prod_{r\leq B_1}r$ and $m_2=\prod_{B_1<r\leq B_2}r$ for primes $r$. The parameters $B_1$ and $B_2$ will be optimized later. For now, we assume that $m_1,m_2<\Lambda$, but $M:=m_1\cdot m_2 \gg \Lambda$, where $\Lambda$ serves the same purpose as in \Cref{alg:sieveimp}. In the sense of the minimization problem (\ref{eq:crem}), the goal of the following algorithm is to compute the elements of the set $\{s\in\bar{\mathcal{L}}_{M}:s<\Lambda\}$. 

\begin{algorithm}
	\emph{Input:} A composite, odd number $N$ and a bound $\Lambda\in\N$
	
	\emph{Output:} A proper factor of $N$ or `$\Lambda$ is too small'
	
	\begin{algorithmic}[1]
		\State Choose bounds $B_1,B_2\in\N$ such that $B_1<B_2$ and $M=m_1\cdot m_2>\Lambda$, where $m_1=\prod_{2<r\leq B_1} r$ and  $m_2=\prod_{B_1<r\leq B_2}r$.
		\State Compute $\gcd(N,M)$. If a factor of $N$ is found, print and stop.
		\State Choose $a,b\in\N$, $a,b<N$ with $\gcd(ab,MN)=1$. 
		\State Let $L:=\lceil 2\sqrt{abN}\rceil$ and $k:= ab\mod M$.
		\For {$r\mid M$} 
		\State Compute $\bar{\mathcal{L}}_{r}:=\{s-L\pmod m: s\in\mathcal{L}_{N,r,k}\}$. 
		\EndFor
		\State Run \Cref{alg:fastcrm} to compute $\bar{\mathcal{L}}_{m_1}$ and $\bar{\mathcal{L}}_{m_2}$ from the sets $\bar{\mathcal{L}}_{r}$.
		\State For all $x_1\in\bar{\mathcal{L}}_{m_1}$, set $j_1:=x_1\cdot m_2^{-1}\pmod{m_1}$. Compute the sorted list $N_1:=\{j_1\cdot m_2\}$. 
		\State Similarly, compute and sort $N_2:=\{j_2\cdot m_1\}$, where $j_2:=x_2\cdot m_1^{-1}\pmod{m_2}$ for $x_2\in\bar{\mathcal{L}}_{m_2}$.
		\For {$\alpha_1 \in N_1$}
		\For {$\alpha_2\in N_2$ with $0\leq \alpha_2<\Lambda - \alpha_1$ or $M-\alpha_1 \leq \alpha_2 <M+\Lambda -\alpha_1$}
		\State Define $\alpha:=\alpha_1+\alpha_2\pmod M$.
		\If{$(\alpha+L)^2-4abN = y^2$ for some $y\in\N$} print $\gcd(\alpha+L-y,N)$.
		\EndIf
		\EndFor
		\EndFor
		\State If no proper factor is found, return `$\Lambda$ is too small'.
	\end{algorithmic}
	\label{alg:tradeoff}
\end{algorithm}

{\prop{		
		If 
		$
		\Lambda>\Lambda_{a,b},
		$
		then \Cref{alg:tradeoff} returns a proper factor of $N$. Ignoring terms that are polynomial in $\log N$ and $\log \Lambda$, the runtime complexity is given by 
		\[
		\widetilde{O}(|N_1|+|N_2|+|\bar{\mathcal{L}}|)
		\]
		bit operations, where $\bar{\mathcal{L}}=\{s\in\mathcal{L}_M: s<\Lambda\}$.
}\label{prop:tradeoff}}

\begin{proof}
	Let us start by proving correctness. The first part of the proof is similar to the one of 
	\Cref{prop:sieveimp}. The choice of parameters in Step 1--4 ensures that $\gcd(abN,M)=1$ and 
	$\gcd(ab,N)=1$, and it is clear that we may compute $\bar{\mathcal{L}}_{m_1}$ and 
	$\bar{\mathcal{L}}_{m_2}$ from the sets $\bar{\mathcal{L}}_{r}$.  Our assumption 
	$\Lambda>\Lambda_{a,b}$ implies that there exists a nontrivial decomposition $N=uv$ of $N$ such 
	that, for $z\in\N$ with $z+L=au+bv$,
	\[
	M>\Lambda>\frac{(au-bv)^2}{4\sqrt{abN}}>z.
	\]
	Note that $z=z\pmod M$ is an element of the set $\{s\in\mathcal{L}_M: s<\Lambda\}$. Moreover, $z\pmod{m_1}\in\bar{\mathcal{L}}_{m_1}$ and $z\pmod{m_2}\in\bar{\mathcal{L}}_{m_2}$. Let $z_1\in N_1$ and $z_2\in N_2$ be the two elements corresponding to these residues in Step 8 and 9. Note that we have $z_1,z_2\leq M$, from which it follows that $z_1+z_2 \leq 2M$. In addition, we know that $z_1+z_2\pmod M < \Lambda$. Clearly, this is satisfied if and only if either $z_1+z_2<\Lambda$ or $M\leq z_1+z_2< M+\Lambda$. These are exactly the intervals we are searching in Step 11, so we will certainly find $z_1+z_2\pmod M$. The Chinese Remainder Theorem implies that 
	\[
	z=z\pmod M =(z_1+z_2)\pmod M.
	\]
	Hence, for the candidate $\alpha=z$, we will retrieve $y=au-bv$ in Step 13 and compute the GCD of $\alpha+L-y=2bv$ and $N$. Since $N$ is odd and $\gcd(ab,N)=1$, this yields a proper factor of $N$. 
	
	Let us now consider the runtime complexity. As in the proof of \Cref{prop:sieveimp}, the runtime complexity for the Steps 1--6 is polynomial in $\log N$ and $\log \Lambda$. Steps 7--9 run in $\widetilde{O}(\max\{|N_1|,|N_2|)$, where we have used that sorting lists of $n$ elements takes $\widetilde{O}(n)$ bit operations (e.g., see \cite[Chap. 2.2]{SedWay}). Let us now consider the for-loop in the Steps 10--13. As we have already explained in the correctness proof, our goal is to find sums of elements of $N_1$ and $N_2$ with $\alpha_1+\alpha_2\pmod M<\Lambda$. For each $\alpha_1\in N_1$, we have to find the values $\alpha_2$ in the intervals in Step 11. However, since $N_2$ is a sorted list, one may easily find the respective area of values $\alpha_2 \in N_2$ that lie in these intervals. We conclude that the Steps 10--13 run in $\widetilde{O}(|N_1|+|\mathcal{L}|)$, which finishes the proof.	
\end{proof}

This extension of the hyperbolic sieve allows to speed-up the technique in \Cref{sec:fermat1}. In the following, we will discuss practical solutions for choosing the parameters of the algorithm in order to minimize the complexity.

\textbf{How to build $m_1$ and $m_2$.} \Cref{prop:tradeoff} implies that we should choose $m_1$ and $m_2$ in a way such that the sets $N_1$, $N_2$ and $\mathcal{L}$ are about the same size. According to \Cref{thm:sievebound}, the size of $N_1$ is equal to
$
|\mathcal{L}_{N,m_1}|\in \widetilde{O}\left(m_1/2^{\omega(m_1)}\right).
$
Repeating the argument of (\ref{eq:extbound}) with a general value of $\Lambda$ (instead of $\lceil \sqrt{N}\rceil$), we assume that the cardinality of $\mathcal{L}$ is bounded by $\widetilde{O}\left(\Lambda/2^{\omega(M)}\right)$. It follows that $N_1$ and $\mathcal{L}$ are about the same size if and only if
$
2^{\omega(M)-\omega(m_1)}\approx\Lambda/m_1,
$
which is equivalent to 
\begin{align}
\omega(m_2)\approx \log_2(\Lambda/m_1).
\label{eq:m2}
\end{align}
For example, we could choose $m_1$ very close to $\Lambda$, forcing $m_2$ to be really small. The extreme case would be $m_1\approx\Lambda$, in which \Cref{alg:tradeoff} would basically degenerate into \Cref{alg:sieveimp}. Or we could choose $m_1$ small, leading to a rather large value for $m_2$. For optimizing the runtime complexity, our goal is to achieve $|N_1|\approx |N_2|$. One way to obtain this in practice is by considering the value $m$ used in \Cref{alg:sieveimp}, for which we have $m\approx \Lambda$. We then define $m_1$ by removing the largest prime factor of $m$, and build $m_2$ according to (\ref{eq:m2}), i.e., by building the product of the $\lceil \log _2(\Lambda/m_1)\rceil$ primes following those contained in $m_1$. By using \Cref{lem:crt}, we may quickly compute the precise cardinalities of the resulting sets $N_1$ and $N_2$. If $|N_2|\ll |N_1|$, we repeat the process by removing further primes from $m_1$ and adjusting $m_2$. We proceed in this manner until we obtain $|N_2|\approx |N_1|$. Also note that we do not necessarily need to build $m_1$ and $m_2$ as products of consecutive primes. We only need to make sure that they are coprime.

{\example{Let $N=17344343992304993085649094809$ like in \Cref{exa:powers}. From the modulus $m=55870214400$ we used there, we remove the prime factor $19$ and define $m_1=2^8\cdot 3^3\cdot 5^2\cdot 7\cdot 11\cdot 13\cdot 17$. Assuming that $\Lambda\approx m$ for our target bound $\Lambda$ and considering (\ref{eq:m2}), we set
\[
\omega(m_2)=\lceil\log_2(m/m_1)\rceil=5.
\]		
We hence define $m_2=19\cdot 23\cdot 29\cdot 31\cdot 37$. For this choice, $N_1$ has $215040$ elements and $N_2$ has $399168$ elements. Running \Cref{alg:tradeoff} returns the factorization of $N$ after checking $138342$ candidates $\alpha$ in Step 13. Comparing this to the $1935360$ elements we had to check in \Cref{exa:powers}, we conclude that \Cref{alg:tradeoff} allows for a notable speed-up of the runtime. In essence, we have traded the factor $19$ from $m_1$ to $m_2$ in return for higher space complexity. We want to stress that the square-root check is much more expensive than building the sets $N_1$ and $N_2$, hence we particularly benefit from the reduction of potential candidates $\alpha$ for $z$.
	}\label{exa:trade}}

{\remark{Due to the heuristic nature of the approach and the dependence on (\ref{eq:extbound}), we are not able to give a rigorous asymptotic runtime analysis of \Cref{alg:tradeoff}. If we were to prove (\ref{eq:extbound}), we might be able to improve upon the constant $C$ of the runtime established in \Cref{thm:main}. In any case, note that the space complexity is non-negligible due to the fact that we have to precompute the sets $N_1$ and $N_2$. A possible approach for reducing the required space is the well-known Schroeppel-Shamir algorithm (see \cite{BecCorJou} and \cite{HowJou}). However, further investigations of potential reductions of the time and space complexity remain a prospect for future research.}
\vspace{12pt}

\textbf{Time-space tradeoffs for semiprime numbers.}
For $N=pq$ with two distinct primes $p$ and $q$, we may further improve the runtime complexity by applying ideas from recent publications on deterministic integer factorization algorithms (\cite{Har},\cite{Hit1},\cite{Hit2}). Let us focus on the case $a=b=1$, which has been considered in \cite{Hit1}. The approach exploits the fact that $N+1\equiv p+q\mod \phi(N)$, which implies
$
\alpha^{p+q}\equiv \alpha^{N+1}\mod N
$
for every $\alpha\in\Z_N^*$. Considering the discussion in the beginning of \Cref{sec:fermat1}, we deduce that our desired value $z<\Lambda$ satisfies
\[
\alpha^z\equiv \alpha^{N+1-L}\mod N,
\]
where $L=\lceil 2\sqrt{N}\rceil$. We may easily compute the value on the right-hand side, and obtain a discrete logarithm problem that can be solved in time $\widetilde{O}(\Lambda^{1/2})$ by applying Shanks' well-known babystep-giantstep method. In \cite{Hit1}, the hyperbolic sieve has first been applied to further improve this runtime by a sub-exponential factor. In addition, the whole approach has been extended to general $a,b\in\N$ in \cite{Hit2}.

While this improvement by taking the square root is quite substantial, it has to be stressed again that it only works for semiprime numbers. For general $N$, we have to use the algorithms presented in this paper.

\section{Appendix: Proofs of the sieve bound modulo prime powers}\label{sec:appendix}
In the following, we prove \Cref{thm:sievepower} on the cardinalities of sieve sets modulo prime powers.

{\definition{Let $N,k,e\in\N$ and $r$ be an odd prime such that $\gcd(Nk,r)=1$. For 
$a\in\mathcal{L}_{N,r^e,k}$, we define
\[
\mathcal{P}_a:=\{(x,y)\in\mathcal{H}_{N,r^e} : a\equiv kx+y \mod r^e\}. 
\]
 For each pair $(x,y)\in \mathcal{P}_a$, we consider the $r$-adic evaluation $\nu_r(kx-y)$, i.e., the maximal exponent of $r$ dividing $kx-y$. We then define
\[
\nu_a:=\min\{\nu_r(kx-y):(x,y)\in \mathcal{P}_a\},
\] 
which we will also call the ``$\nu$-value of $a$''.
}\label{def:va}}
\vspace{12pt}

The following lemma concerns the relations between elements of the sieve sets modulo prime powers and thus prepares the proof for \Cref{thm:sievepower}.

{\lemma{Let $N,k,e\in\N$ and $r$ be an odd prime such that $\gcd(Nk,r)=1$. Then each element in $\mathcal{L}_{N,r^{e+1},k}$ reduces to an element in $\mathcal{L}_{N,r^{e},k}$ modulo $r^e$. For $a\in \mathcal{L}_{N,r^{e},k}$, let $\#_a$ denote the number of elements in  $\mathcal{L}_{N,r^{e+1},k}$ that reduce to $a$ modulo $r^e$. Then
		\[
			\#_a=
		\begin{cases}
			r     & \mbox{if } \nu_a<e/2,\\
			(r+1)/2    & \mbox{if } \nu_a=e/2,\\
			1   & \mbox{if } \nu_a>e/2.
		\end{cases}
		\]
}\label{lem:acount}}
\begin{proof}
	Let $b\in\mathcal{L}_{N,r^{e+1},k}$ be arbitrary. By definition, there exists $(x,y)$ in 
	$\mathcal{H}_{N,r^{e+1}}$ such that $b=kx+y\pmod{r^{e+1}}$. It is also clear that there are 
	$s,t\in\{0,1,\ldots,r-1\}$ such that $x=\bar{x}+sr^e$ and $y=\bar{y}+tr^e$ for some 
	$(\bar{x},\bar{y})\in\mathcal{H}_{N,r^{e}}$. Hence, we get
	\[
	kx+y=k(\bar{x}+sr^e)+\bar{y}+tr^e=k\bar{x}+\bar{y}+r^e(ks+t),
	\]
	and it follows that $b\equiv kx+y\equiv k\bar{x}+\bar{y}\mod r^e$. The element $k\bar{x}+\bar{y}\pmod{r^e}$ is in  $\mathcal{L}_{N,r^{e},k}$, so we have shown the first claim.
	
 	We now consider $a\in\mathcal{L}_{N,r^{e},k}$. We set $b_j:=a+jr^e$ for $j\in\{0,1,\ldots,r-1\}$ and count the number of values $j$ for which $b_j$ is in $\mathcal{L}_{N,r^{e+1},k}$. We already know that $b_j$ needs to be of the form $kx+y$ for $x=\bar{x}+sr^e$ and $y=\bar{y}+tr^e$, where $(\bar{x},\bar{y})$ is some pair in $\mathcal{H}_{N,r^{e}}$ such that $a\equiv k\bar{x}+\bar{y} \mod r^e$. Moreover, we know that we have
	\begin{equation}\label{eq:1a}
		j\equiv ks+t \mod r.
	\end{equation}
	In addition, we require $xy\equiv N \mod r^{e+1}$. For $c:=(\bar{x}\bar{y}-N)/r^e$, we obtain
	\[
	N\equiv xy\equiv (\bar{x}+sr^e)(\bar{y}+tr^e)\equiv \bar{x}\bar{y}+r^e(s\bar{y}+t\bar{x})\equiv N+r^e(c+s\bar{y}+t\bar{x})\mod r^{e+1}.
	\]
	Clearly, this holds if and only if 
	\begin{equation}\label{eq:2a}
		-c\equiv s\bar{y}+t\bar{x}\mod r.
	\end{equation}
	We consider $(\ref{eq:1a})$ and $(\ref{eq:2a})$ as two linear equations in the variables $s$ and $t$. By some simple transformations,  we can reduce them to $\bar{x} j - c\equiv s(k\bar{x}-\bar{y})\mod r$. One easily observes that there exists a solution $s$ if and only if 
	\begin{equation}\label{eq:3}
	k\bar{x}\not\equiv \bar{y}\mod r \,\,\, \vee \,\,\, j\equiv c\bar{x}^{-1}\mod r.
	\end{equation}
If the first condition is satisfied, the equations $(\ref{eq:1a})$ and $(\ref{eq:2a})$ are always solvable (for each choice of $j$), and $b_j$ is an element of $\mathcal{L}_{N,r^{e+1},k}$ for all $j\in\{0,1,\ldots,r-1\}$. If the first condition is not satisfied, then the only value for which we have a solution is $j\equiv c\bar{x}^{-1}\mod r$. In this case, we need to count the different possibilities to represent $a\equiv k\bar{x}+\bar{y} \mod r^e$ with pairs $(\bar{x},\bar{y})\in\mathcal{H}_{N,r^{e}}$, as each of these pairs may give rise to a different value of $j\in\{0,1,\ldots,r-1\}$. In \Cref{def:va}, we have denoted the set of these pairs as $\mathcal{P}_a$.

Let us start with the case $\nu_a=0$. This implies that $k\bar{x}\not\equiv \bar{y}\mod r$ for every pair $(\bar{x},\bar{y})\in\mathcal{P}_a$. As we have discussed above, $b_j$ is an element of $\mathcal{L}_{N,r^{e+1},k}$ for every $j\in\{0,1,\ldots,r-1\}$, so $a\in\mathcal{L}_{N,r^e,k}$ induces $r$ elements in $\mathcal{L}_{N,r^{e+1},k}$.

Assume now that $\nu_a>0$. Then the first condition in (\ref{eq:3}) is not true. As explained, the resulting number of possible values for $j$ depends on the number of different possibilities to represent $a$, i.e., the cardinality of $\mathcal{P}_a$. It is clear that $k\bar{x}\equiv \bar{y}\mod r^{\nu_a}$ for every $(\bar{x},\bar{y})\in \mathcal{P}_a$. For any two pairs $(x_0,y_0),(x_1,y_1)\in\mathcal{P}_a$, it follows that
\[
2y_0\equiv kx_0+y_0\equiv a \equiv kx_1+y_1\equiv 2y_1 \mod r^{\nu_a},
\]
hence $y_0\equiv y_1\mod r^{\nu_a}$. This implies that all pairs in $\mathcal{P}_a$ are equal modulo $r^{\nu_a}$. In particular, this is true for the values of $\bar{x}^{-1}\pmod r$. So in order to count the number of different values of $j\equiv c\bar{x}^{-1}\mod r$ we can get, it suffices to count the different values of $c$ related to the pairs in $\mathcal{P}_a$. We proceed by considering different sizes of $\nu_a$.

We start with $\nu_a>e/2$. Let $(\bar{x},\bar{y})\in \mathcal{P}_a$ be arbitrary. We have $N\equiv \bar{x}\bar{y}\mod r^e$ and want to extend $(\bar{x},\bar{y})$ to the other pairs in $\mathcal{P}_a$. We have discussed above that all these pairs are equal modulo $r^{\nu_a}$. We may hence write them as
\[
(\bar{x}-ir^{\nu_a})(\bar{y}+kir^{\nu_a})\equiv N+ir^{\nu_a}(k\bar{x}-\bar{y})-ki^2r^{2\nu_a}\equiv N  \mod r^e.
\]
We thus have $(\bar{x}-ir^{\nu_a},\bar{y}+kir^{\nu_a})\in\mathcal{H}_{N,r^e}$, and one easily observes that their corresponding element in $\mathcal{L}_{N,r^e,k}$ is indeed equal to $a$.  Let us now consider the values for $c$ that correspond to these pairs. For $\gamma=r^{\nu_a-e}(k\bar{x}-\bar{y})$ and $\delta=(\bar{x}\bar{y}-N)/r^{e}$, we obtain
\[
c_i\equiv \frac{(\bar{x}-ir^{\nu_a})(\bar{y}+kir^{\nu_a})-N}{r^{e}}\equiv \delta+\gamma i-r^{2\nu_a-e}ki^2\equiv \delta\mod r,
\]
since $\gamma$ is still divisible by $r$. So the value of $c$ is actually equal to $\delta$ for all pairs in $\mathcal{P}_a$. As a result, the elements $a\in\mathcal{L}_{N,r^e,k}$ with $\nu_a>e/2$ only induce one element in $\mathcal{L}_{N,r^{e+1},k}$. 

Next, we consider the case $1\leq\nu_a\leq e/2$. Again, let $(\bar{x},\bar{y})\in \mathcal{P}_a$ be arbitrary. The only way we can extend the pair $(\bar{x},\bar{y})$ to other pairs in $\mathcal{H}_{N,r^e}$ that are equal modulo $r^{\nu_a}$ is by extending them modulo $r^{e-\nu_a}$, i.e.
\[
(\bar{x}-ir^{e-\nu_a})(\bar{y}+kir^{e-\nu_a})\equiv N+ir^{e-\nu_a}(k\bar{x}-\bar{y})-ki^2r^{2(e-\nu_a)}\equiv N  \mod r^e,
\]
and, for $\gamma=r^{-\nu_a}(k\bar{x}-\bar{y})$ and $\delta$ as above,
\begin{equation}
c_i\equiv \frac{(\bar{x}-ir^{e-\nu_a})(\bar{y}+kir^{e-\nu_a})-N}{r^{e}}\equiv \delta+\gamma i-r^{e-2\nu_a}ki^2\mod r.
\label{eq:4}
\end{equation}
If $\nu_a<e/2$, then $c_i\equiv \delta+\gamma i \mod r$ and we have $r$ possible distinct values for $c_i$, as $\gamma$ is not divisible by $r$. Therefore, the elements $a\in\mathcal{L}_{N,r^e,k}$ with $\nu_a<e/2$ induce $r$ elements in $\mathcal{L}_{N,r^{e+1},k}$. We are left with $\nu_a=e/2$. In this case, (\ref{eq:4}) is a quadratic equation in $i$ with the discriminant 
\[
\gamma^2+4k(\delta-c_i).
\] 
It is solvable whenever $c_i$ takes on a value such that the discriminant is a quadratic residue (including $0$). We conclude that there are $(r+1)/2$ distinct possible values for $c_i$, and that the $a\in\mathcal{L}_{N,r^e,k}$ with $\nu_a=e/2$ hence induce $(r+1)/2$ elements in $\mathcal{L}_{N,r^{e+1},k}$.
\end{proof}

{\lemma{Let $N,k,e\in\N$ and $r$ be an odd prime such that $\gcd(Nk,r)=1$. In addition, assume that $(k^{-1}N|r)=1$. Then there exist exactly two elements $a,b\in\mathcal{L}_{N,r^e,k}$ with $\nu_a,\nu_b\geq e/2$, both satisfying $\nu_a=\nu_b=\lceil e/2\rceil$.}\label{lem:powcount}}

\begin{proof}	
	We prove this by induction on $e$. We consider three base cases and start with $e=1$. There are 
	two distinct solutions of 
	\begin{equation}\label{eq:6}
	x^2\equiv k^{-1}N\mod r,
	\end{equation}
	both satisfying $x^2\not\equiv k^{-1}N\mod r^2$. One easily checks that they correspond to two distinct elements $a,b\in\mathcal{L}_{N,r,k}$ satisfy $\nu_a=\nu_b=1=\lceil e/2\rceil$.
	
	We now consider $e=2$. Besides the two elements $a,b\in\mathcal{L}_{N,r,k}$ with $\nu_a=\nu_b=1$, all other elements in this set have a $\nu$-value of $0$. Let $c$ be such an element with $\nu_c=0$, and let $s\in\mathcal{L}_{N,r^2,k}$ be induced by $c$. Then each $(x_s,y_s)\in\mathcal{P}_s$ reduces to some $(x_c,y_c)\in\mathcal{P}_c$ modulo $r$, and the latter pairs all satisfy $x_c^2\not\equiv k^{-1}N\mod r$. Hence, we have shown that $\nu_s=0$. We may thus focus on the elements induced by $a$ and $b$. \Cref{lem:acount} implies that they both induce exactly one element in $\mathcal{L}_{N,r^2,k}$. Let $a^*$ be the element that is induced by $a$. In addition, let $x_a$ be the solution of (\ref{eq:6}) corresponding to $a$, and $y_a=x_a^{-1}N\mod r$. For $y^*=x_a^{-1}N\mod r^2$, we obtain $y^*\equiv y_a\mod r$, and it follows that the pair $(x_a,y^*)$ is in $\mathcal{P}_{a^*}$. Since we know that $x_a^2\not\equiv k^{-1}N\mod r^2$, we conclude that $\nu_{a^*}=1=\lceil e/2\rceil$. Repeating this argument for $b$ and its induced element $b^*$, we also obtain $\nu_{b^*}=1=\lceil e/2\rceil$.
	
	We consider the final base case $e=3$. Again, we may focus on the two elements $a,b\in\mathcal{L}_{N,r^2,k}$ with $\nu_a=\nu_b=1$. However, \Cref{lem:acount} now implies that they both induce $(r+1)/2$ elements in $\mathcal{L}_{N,r^3,k}$. Let us focus on $a$, and let $x_a$ be the corresponding solution of (\ref{eq:6}) as above. Then one observes that there is a unique value $i_0\in\{0,1,\ldots,r-1\}$ such that
	\[
	(x_a+i_0r)^2\equiv k^{-1} N\mod r^{2}.
	\] 
	It is also clear that one of the pairs in $\mathcal{P}_a$ is of the form $(x_a+i_0r,y^*)$, where $y^*=(x_a+i_0r)^{-1}N\pmod{r^2}$. As we have seen in the proof of \Cref{lem:acount}, there is exactly one element $a^*$ of the $(r+1)/2$ elements in $\mathcal{L}_{N,r^3,k}$ induced by $a$ that corresponds to this particular pair. Therefore, each pair $(x,y)\in\mathcal{P}_{a^*}$ satisfies
	$
	x^2\equiv (x_a+i_0r)^2\equiv k^{-1} N\mod r^2.
	$
	Also, we clearly have $(x_a+i_0r,y^*)\in\mathcal{P}_{a^*}$ for $y^*=(x_a+i_0r)^{-1}N\pmod{r^3}$, and one observes that
	$(x_a+i_0r)^2\not\equiv k^{-1} N\mod r^{3}$. We deduce $\nu_{a^*}=2=\lceil e/2\rceil$.  All the other elements $s_i$ induced by $a$ correspond to other pairs in $\mathcal{P}_a$, for which the first component is of the form $x_a+ir$ with $i\neq i_0$.  Hence, for each $(x,y)\in\mathcal{P}_{s_i}$, we have $x^2\equiv (x_a+ir)^2\not\equiv k^{-1} N\mod r^2$. It follows that the $\nu$-value of these other elements induced by $a$ is still equal to $1<e/2$. Repeating this argument for $b$ concludes this base case.
	
	We may now describe the induction step for the general case. First, assume that $e=2n$ for $n\in\N$. According to the induction assumption, there are exactly two elements $a,b\in\mathcal{L}_{N,r^{2n-1},k}$ with $\nu$-value equal to $\lceil (e-1)/2\rceil=n$. We note that there are two distinct solutions $x_a$ and $x_b$ of 
	$
	x^2\equiv k^{-1}N\mod r^n
	$
	which correspond to $a$ and $b$, in the sense that $(x_a,y_a^*)\in\mathcal{P}_a$ and $(x_b,y_b^*)\in\mathcal{P}_b$ for suitable $y_a^*$ and $y_b^*$. By similar arguments as in the base case, we easily observe that the $\nu$-value of those elements in $\mathcal{L}_{N,r^{2n-1},k}$ that are not induced by $a$ and $b$ is still strictly smaller than $n$, and thus smaller or equal to $n-1<e/2$. \Cref{lem:acount} implies that $a$ induces exactly one element $a^*$ in $\mathcal{L}_{N,r^{2n},k}$. Moreover, $x_a^2\not \equiv k^{-1}N\mod r^{n+1}$ implies
	\[
	\nu_{a^*}=\lceil (e-1)/2\rceil=n= \lceil e/2\rceil.
	\]
	The same is true for the element $b^*$ induced by $b$, so we may conclude the proof for $e=2n$.	Finally, let us consider $e=2n-1$ for $n\in\N$. Again, the induction assumption implies that there are exactly two elements $a,b\in\mathcal{L}_{N,r^{2n-2},k}$ with $\nu$-value equal to $\lceil (e-1)/2\rceil=n-1$. \Cref{lem:acount} implies that both induce $(r+1)/2$ elements in $\mathcal{L}_{N,r^{2n-1},k}$. In addition, there is a unique value $i_0\in\{0,1,\ldots,r-1\}$ with 
	\[
	(x_a+i_0r^n)^2\equiv k^{-1}N\mod r^{n+1},
	\] 
	and a similar statement holds for $x_b$. By retracing the argument in the base case $e=3$, we obtain $a^*,b^*\in\mathcal{L}_{N,r^{2n-1},k}$ with 
	\[
	\nu_{a^*}=\nu_{b^*}=\lceil (e-1)/2\rceil+1=n-1+1=\lceil e/2\rceil.
	\]
	The argument also shows that the $\nu$-value of the other elements is still equal to or smaller than $\lceil (e-1)/2\rceil=n-1<e/2$. This concludes the proof.
\end{proof}

\begin{proof}[Proof of \Cref{thm:sievepower}]	
	We start with the case $(k^{-1}N|r)=-1$. Let $a\in\mathcal{L}_{N,r^e,k}$ be arbitrary. For every $(\bar{x},\bar{y})\in\mathcal{P}_a$, we deduce that $k\bar{x}\not\equiv \bar{y}\mod r$, as otherwise we would have $\bar{x}^2\equiv k^{-1}N\mod r$. Hence, $\nu_a=0$ for every $a\in \mathcal{L}_{N,r^e,k}$, and \Cref{lem:acount} implies that the first claim of \Cref{thm:sievepower} is correct.
	
	Let us proceed with $(k^{-1}N|r)=1$. We consider even and odd values of $e$ separately. Let us start with $e=2n-1$ for $n\in\N$. For any $a\in\mathcal{L}_{N,r^e,k}$, we either have $\nu_a<e/2$ or $\nu_a>e/2$. \Cref{lem:acount} implies that those elements with $\nu_a<e/2$ each induce $r$ elements in $\mathcal{L}_{N,r^{e+1},k}$, while those with $\nu_a>e/2$ each induce only one element in $\mathcal{L}_{N,r^{e+1},k}$. \Cref{lem:powcount} implies that there are exactly two such elements $a,b$ with $\nu_a=\nu_b>e/2$, so we have to subtract $2(r-1)$ from $|\mathcal{L}_{N,r^e,k}|\cdot r$ to obtain the cardinality of $\mathcal{L}_{N,r^{e+1},k}$.
	
	Finally, let $e=2n$. For any $a\in\mathcal{L}_{N,r^e,k}$, we can have $\nu_a<e/2$ or $\nu_a\geq e/2$. The elements corresponding to the first case each induce $r$ elements in $\mathcal{L}_{N,r^{e+1},k}$. Moreover, \Cref{lem:powcount} implies that there are exactly two elements $a,b\in\mathcal{L}_{N,r^e,k}$ corresponding to the second case, and both satisfy $\nu_a=\nu_b=e/2$. According to \Cref{lem:acount}, they both induce $(r+1)/2$ elements in $\mathcal{L}_{N,r^{e+1},k}$. So for the total number of elements in $\mathcal{L}_{N,r^{e+1},k}$, we obtain
	\[
	|\mathcal{L}_{N,r^e,k}|\cdot r-2(r-1)/2 = |\mathcal{L}_{N,r^e,k}|\cdot r - (r-1).
	\]
\end{proof}

The following statement about shifted quadratic residues modulo prime powers is an immediate consequence of our results. Since we have not found it in the literature, we state it here.

{\corollary{Let $t,e\in\N$ and $r$ be an odd prime such that $\gcd(t,r)=1$. For the set $Q_e$  of quadratic residues modulo $r^e$, we define
$
\mathcal{S}_e:=\{x\in\Z_{r^e}:x^2-t\in Q_e\}.
$
For $\bar{t}:=4^{-1}t\pmod{r^e}$, we have $|\mathcal{S}_1|=(r+(\bar{t}|r))/2$ and
\[
|\mathcal{S}_{e+1}|=
\begin{cases}
	|\mathcal{S}_e|\cdot r     & \mbox{if } (\bar{t}|r)=-1,\\
	|\mathcal{S}_e|\cdot r-2^{e\pmod 2}\cdot (r-1)     & \mbox{if } (\bar{t}|r)=1.
\end{cases}
\]
}}
\begin{proof}
\Cref{lem:qres} implies that $\mathcal{S}_e=\mathcal{L}_{\bar{t},r^e,1}$. The statement then follows from \Cref{lem:card} and \Cref{thm:sievepower}.
\end{proof}

{\remark{Let $N,k,e\in\N$ with $\gcd(Nk,2)=1$. According to our computations, the cardinality of the sieve sets modulo powers of $2$ is as follows. We have
		$|\mathcal{L}_{N,2,k}|=|\mathcal{L}_{N,4,k}|=1$ and
		\[
		|\mathcal{L}_{N,8,k}|=
		\begin{cases}
			1     & \mbox{if } N\equiv 3k \mod 4,\\
			2    & \mbox{else}.\\
		\end{cases}
		\]
		Moreover, we have $|\mathcal{L}_{N,16,k}|=2$ and
		\[
		|\mathcal{L}_{N,32,k}|=
		\begin{cases}
			2    & \mbox{if } N\equiv 5k\mod 8,\\
			4   & \mbox{else}.
		\end{cases}
		\]
		For $e\geq 5$, we have
		\[
		|\mathcal{L}_{N,2^{e+1},k}|=
		\begin{cases}
			|\mathcal{L}_{N,2^e,k}|\cdot 2     & \mbox{if } N\equiv 3k \mod 4,\\
			|\mathcal{L}_{N,2^e,k}|\cdot 2    & \mbox{if } N\equiv 5k\mod 8,\\
			\left(|\mathcal{L}_{N,2^e,k}| -2^{e\pmod 2}\right)\cdot 2   & \mbox{if } N\equiv k \mod 
			8.
		\end{cases}
		\]
Most likely, this can be proved via slight modifications to the strategy for proving the general result for odd prime powers.
}\label{rem:twopower}}

{\remark{Assume that we know that $u^2\equiv k^{-1} N\mod r^n$ or, equivalently,
		\[
		ku\equiv v \mod r^n
		\]
		holds for some odd prime $r$, some $n\in\N$ and some co-divisors $u,v$ of $N$. Then one easily deduces from \Cref{lem:powcount} that the residue of $ku+v$ modulo $r^{2n-1}$ must be one of the two values in $a\in\mathcal{L}_{N,r^{2n-1},k}$ with $\nu_a=n$. As a consequence, knowing this information about the divisors of $N$ modulo $r^n$ allows to deduce two candidates for the corresponding linear combination $ku+v$ modulo $r^{2n-1}$. The two possible values $a$ may be easily computed by following the arguments in the proof of \Cref{lem:powcount}.
		While this fact is certainly interesting, it is clear that it strongly depends on the 
		congruence $ku\equiv v\mod r^n$ for our choice of the value of $k$. For instance, we can 
		apply this idea in the search for $u+v$ only if we know any odd prime powers dividing the 
		divisor difference $u-v$. However, the approach certainly allows for improvements in cases 
		where the divisors are known to be of certain shape. For instance, this is the case for 
		(generalized) Mersenne and Fermat numbers.
}}


\end{document}